\documentclass[useAMS]{gITR2e}

\newcommand{\bfx}{{\bf x}}
\newcommand{\bfxp}{{{\bf x}^\prime}}
\newcommand{\N}{{\mathbf N}}

\newcommand{\Q}{{\mathbf Q}}
\newcommand{\R}{{\mathbf R}}

\newcommand{\g}{{\mathfrak g}}
\newcommand{\li}{{\mathfrak l}}
\newcommand{\hii}{{\mathfrak h}}

\newcommand{\Z}{{\mathbf Z}}
\newcommand{\C}{{\mathbf C}}

\newcommand{\mcg}{{\mathcal G}}

\def\cprime{$'$}

\makeatletter
\def\eqnarray{\stepcounter{equation}\let\@currentlabel=\theequation
\global\@eqnswtrue
\tabskip\@centering\let\\=\@eqncr
$$\halign to \displaywidth\bgroup\hfil\global\@eqcnt\z@
  $\displaystyle\tabskip\z@{##}$&\global\@eqcnt\@ne
  \hfil$\displaystyle{{}##{}}$\hfil
  &\global\@eqcnt\tw@ $\displaystyle{##}$\hfil
  \tabskip\@centering&\llap{##}\tabskip\z@\cr}

\def\endeqnarray{\@@eqncr\egroup
      \global\advance\c@equation\m@ne$$\global\@ignoretrue}

\def\@yeqncr{\@ifnextchar [{\@xeqncr}{\@xeqncr[5pt]}}
\makeatother

\begin{document}
\doi{10.1080/10652460YYxxxxxxx}
 \issn{1476-8291}
\issnp{1065-2469}
 \jvol{00} \jnum{00} \jyear{2009} \jmonth{January}

\markboth{Howard S.~Cohl}{Integral Transforms and Special Functions}


\title{
Fourier expansions for a logarithmic fundamental
solution of the polyharmonic equation
}

\author{Howard S.~Cohl$^{\rm a,b}$$^{\ast}$\thanks{$^\ast$Corresponding author. Email: howard.cohl@nist.gov
\vspace{6pt}}
\\\vspace{6pt}  $^{\rm a}${\em{Applied and Computational Mathematics Division,
Information Technology Laboratory,
National Institute of Standards and Technology,
Gaithersburg, Maryland, U.S.A.
}}; \\$^{\rm b}${\em{Department of Mathematics, University of Auckland, Auckland, New Zealand
}}\\\vspace{6pt}\received{ v3.5 released October 2008} }

\maketitle

\begin{abstract}
In even-dimensional Euclidean space for integer powers of the Laplacian greater than
or equal to the dimension divided by two, a fundamental solution for the polyharmonic
equation has logarithmic behavior.  We give two approaches for developing a Fourier
expansion of this logarithmic fundamental solution. The first approach is algebraic
and relies upon the construction of two-parameter polynomials.  We
describe some of the properties of these polynomials, and use them to derive 
the Fourier expansion for a logarithmic fundamental solution of the polyharmonic equation.
The second approach depends on the computation of parameter derivatives of
Fourier series for a power-law fundamental solution of the polyharmonic equation. 
The resulting Fourier series is given in terms of sums over associated
Legendre functions of the first kind.  We conclude by comparing the two approaches and
giving the azimuthal Fourier series for a logarithmic fundamental solution of
the polyharmonic equation in rotationally-invariant coordinate systems.\bigskip

\begin{keywords}
fundamental solutions; polyharmonic equation; Fourier series; polynomials; associated
Legendre functions
\end{keywords}
\begin{classcode}
35A08; 31B30; 31C12; 33C05; 42A16
\end{classcode}\bigskip

\end{abstract}

\section{Introduction}
\label{Introduction}

Solutions of the polyharmonic equation (powers of the Laplacian operator)
are ubiquitous in many areas of computational, pure, applied mathematics,
physics and engineering.
We concern ourselves, in this paper, with a fundamental solution of the 
polyharmonic equation (Laplace, biharmonic, etc.), which by convolution 
yields a solution to the 
inhomogeneous polyharmonic equation.  
Solutions to inhomogeneous polyharmonic equations are useful in many 
physical applications including those areas related to Poisson's equation
such as Newtonian gravity, electrostatics, magnetostatics, quantum
direct and exchange interactions
(cf.~\S 1 in Cohl \& Dominici (2010) \cite{CohlDominici}), etc.
Furthermore, applications of higher-powers of the Laplacian
include such varied areas as minimal surfaces \cite{MonterdeUgail}, Continuum 
Mechanics \cite{LaiRubinKrempl},
Mesh deformation \cite{Helenbrook}, Elasticity \cite{LurieVasiliev},
Stokes Flow \cite{Kirby}, Geometric Design \cite{Ugail},
Cubature formulae \cite{Sobolev}, mean value theorems 
(cf.~Pizzetti's formula) \cite{Nicolescu}, and Hartree-Fock calculations 
of nuclei \cite{Vautherin}.  

It is a well-known fact (see for 
instance Schwartz (1950) (\cite{Schw}, p.~45),
Gel{\cprime}fand \& Shilov (1964) (\cite{GelfandShilov}, p.~202)
that a fundamental solution of the polyharmonic equation
on $d$-dimensional Euclidean space is given by combinations of 
power-law and logarithmic functions of the global distance between two points.
In a recent paper (Cohl \& Dominici (2010) \cite{CohlDominici}),
we derived a complex identity which determined
the Fourier coefficients of a power-law fundamental solution
of the polyharmonic equation.  The Fourier coefficients were
seen to be given in terms of associated Legendre functions.
The present work is concerned with computing the Fourier 
coefficients of a logarithmic fundamental solution of the 
polyharmonic equation.  One obtains a logarithmic fundamental 
solution for the polyharmonic equation only in even-dimensional
Euclidean space and only when the power of the Laplacian 
is greater than or equal to the dimension divided by two. 
The most familiar example of a logarithmic fundamental solution of
the polyharmonic equation occurs in two-dimensions, for a 
single-power of the Laplacian, i.e., Laplace's equation.

We present two different approaches for obtaining Fourier series
of a logarithmic fundamental solution for the polyharmonic 
equation. The first approach is algebraic and involves
the generation of a certain set of naturally arising two-index polynomials
which we refer to as logarithmic polynomials.
The second approach starts with the main result from 
Cohl \& Dominici (2010) \cite{CohlDominici} and determines
the Fourier series expansion for a logarithmic fundamental solution
of the polyharmonic equation through parameter differentiation.
Series expansions for fundamental solutions of 
linear partial differential equations such as the polyharmonic equation
are extremely useful in determining 
Dirichlet boundary values for solutions on interior domains
(see for example Cohl \& Tohline (1999) \cite{CT}).

This paper is organized as follows.  
In \S\ref{Prelude} we introduce the problem.
In \S\ref{AlgebraicapproachtologarithmicFourierseries}
we describe our algebraic approach to computing a Fourier
series of a logarithmic fundamental solution of the 
polyharmonic equation.
In \S\ref{LimitderivativeapproachtologarithmicFourierseries} we
give our limit derivative approach for computing the Fourier
series of a logarithmic fundamental solution of the polyharmonic
equation.
In \S\ref{Comparison} we give some comparisons between the two approaches.
In \S\ref{FourierExpansionforlogfundsolofpolyeq} we use the results
presented in the previous sections to obtain azimuthal Fourier expansions for 
a logarithmic fundamental solution of the polyharmonic equation
in rotationally-invariant coordinate systems which parametrize
points in $d$-dimensional Euclidean space.
In Appendix \ref{DerivativeswithrespecttothedegreeofP} we 
present some necessary formulae relating to differentiation of
associated Legendre functions of the first kind with respect to
the degree.
In Appendix \ref{Thelogarithmicpolynomials} we present some of the properties
of the logarithmic polynomials.

Throughout this paper we rely on the following definitions.  
The set of natural numbers is given by $\N:=\{1,2,3,\ldots\}$, the set
$\N_0:=\{0,1,2,\ldots\}=\N\cup\{0\}$, the set of integers
is given by $\Z:=\{0,\pm 1,\pm 2,\ldots\},$ and the set 
$\Q$ represents the rational numbers. 
For $a_1,a_2,\ldots\in\C$, if $i,j\in\Z$ and $j<i$ then
$\sum_{n=i}^{j}a_n=0$ and $\prod_{n=i}^ja_n=1$,
where $\C$ represents the complex numbers.  The set $\R$ 
represents the real numbers.
\section{Fundamental solution of the polyharmonic equation and the non-logarithmic Fourier series}
\label{Prelude}

If $\Phi$ satisfies the polyharmonic equation given by
\begin{equation}
(-\Delta)^k\Phi(\bfx)=0,
\label{polyharmoniceq}
\end{equation}
where $\bfx\in\R^d,$ $\Delta:C^p(\R^d)\to C^{p-2}(\R^d)$, for $p\ge 2$ is the 
Laplacian operator defined by
$\Delta:=\frac{\partial^2}{\partial x_1^2}+\ldots+\frac{\partial^2}{\partial x_d^2},$
$k\in\N$ and $\Phi\in C^{2k}(\R^d),$ then $\Phi$ is called 
polyharmonic. If the power
$k$ of the Laplacian equals two, then (\ref{polyharmoniceq}) is called the 
biharmonic equation and $\Phi$ is called biharmonic.
The inhomogeneous polyharmonic equation is given by
\begin{equation}
(-\Delta)^k\Phi({\bf x})=\rho({\bf x}),
\label{polyh}
\end{equation}
where we take $\rho$ to be an integrable function so that 
a solution to (\ref{polyh}) exists.  A fundamental solution for the 
polyharmonic equation on $\R^d$ is a function 
${\g}_k^d:(\R^d\times\R^d)\setminus\{(\bfx,\bfx):\bfx\in\R^d\}\to\R$ 
which satisfies the equation
\begin{equation}
(-\Delta)^k{\g}_k^d({\bf x},{\bf x}^\prime)=c\delta({\bf x}-{\bf x}^\prime),
\label{unnormalizedfundsolnpolydefn}
\end{equation}
for some $c\in\R, c\ne 0$, where $\delta$ is the Dirac delta function
and $\bfxp\in\R^d$.
When $c=1$, we call a fundamental solution of the polyharmonic equation normalized, 
and denote it by
$\mcg_k^d:(\R^d\times\R^d)\setminus\{(\bfx,\bfx):\bfx\in\R^d\}\to\R$.
The Euclidean inner product $(\cdot,\cdot):\R^d\times\R^d\to\R$ 
defined by
$(\bfx,\bfxp):=x_1x_1^\prime+\ldots+x_dx_d^\prime,$
induces a norm (the Euclidean norm) $\|\cdot\|:\R^d\to[0,\infty)$,
on the finite-dimensional vector space $\R^d$, 
given by $\|\bfx\|:=\sqrt{(\bfx,\bfx)}.$
In the rest of this paper, we will use the gamma function 
$\Gamma:\C\setminus-\N_0\to\C$, which is a natural generalization 
of the factorial function
(see for instance Chapter 5 in Olver {\it et al.}~(2010)
\cite{NIST}).  A fundamental solution of the polyharmonic equation
is given by the following theorem.

\begin{theorem} Let $d,k\in\N$.  Define
\[
\mcg_k^d({\bf x},{\bf x}^\prime)=
\left\{ \begin{array}{ll}

{\displaystyle \frac{(-1)^{k+d/2+1}\ \|{\bf x}-{\bf x}^\prime\|^{2k-d}}
{(k-1)!\ \left(k-d/2\right)!\ 2^{2k-1}\pi^{d/2}}
\left(\log\|{\bf x}-{\bf x}^\prime\|-\beta_{k-d/2,d}\right)}\\[2pt]
\hspace{7.4cm} \mathrm{if}\  d\,\,\mathrm{even},\ k\ge d/2,\\[10pt]
{\displaystyle \frac{\Gamma(d/2-k)\|{\bf x}-{\bf x}^\prime\|^{2k-d}}
{(k-1)!\ 2^{2k}\pi^{d/2}}} \hspace{3.37cm} \mathrm{otherwise},
\end{array} \right. 
\]
where $\beta_{p,d}\in\Q$ is defined as 
$\beta_{p,d}:=\frac12\left[H_p+H_{d/2+p-1}-H_{d/2-1} \right],$
with $H_j$ being the $j$th harmonic number 
\[
H_j:=\sum_{i=1}^j\frac1i,
\]
then $\mcg_k^d$ is a normalized fundamental solution for $(-\Delta)^k$
on Euclidean space $\R^d$.  

\label{greenpoly}
\end{theorem}
\noindent Proof.  See Cohl (2010) \cite{CohlthesisII} and Boyling (1996) \cite{Boyl}.
\medskip

A separable rotationally-invariant coordinate system for the polyharmonic equation
(\ref{polyharmoniceq}) on $\R^d$ is given by
\begin{equation}
\left.
\begin{array}{rcl}
x_1&=&R(\xi_1,\ldots,\xi_{d-1})\cos\phi\\[0.1cm]
x_2&=&R(\xi_1,\ldots,\xi_{d-1})\sin\phi\\[0.1cm]
x_3&=&x_3(\xi_1,\ldots,\xi_{d-1})\\[0.1cm]
&\vdots&\\[0.1cm]
x_d&=&x_d(\xi_1,\ldots,\xi_{d-1})
\end{array}
\quad\right\},
\label{rotatioanallyinvariant}
\end{equation}
which is described by an angle $\phi\in\R$ and $(d-1)$-curvilinear 
coordinates $(\xi_1,\ldots,\xi_{d-1})$. A separable rotationally-invariant coordinate
system transforms the polyharmonic equation into a set of $d$-uncoupled ordinary 
differential equations with separation 
constants $m\in\Z$ and $k_j$ for $1\le j\le d-2$.  For a separable rotationally-invariant 
coordinate system, this uncoupling
is accomplished, in general, by assuming a solution to (\ref{polyharmoniceq}) of the form
\[
\Phi(x)=e^{im\phi}\,{\mathcal R}(\xi_1,\ldots,\xi_{d-1})\prod_{i=1}^{d-1} 
A_i(\xi_i,m,k_1,\ldots,k_{d-2}),
\]
where the domains of the functions 
${\mathcal R}$ and $A_i$, for $1\le i\le d-1$, and the constants $k_j$ for 
$1\le j\le d-1$, depend on the specific rotationally-invariant coordinate system.
A rotationally-invariant coordinate system parametrizes points on the $(d-1)$-dimensional
half-hyperplane given by $\phi=const.$ and $R\ge 0$ using curvilinear coordinates
$(\xi_1,\ldots,\xi_{d-1})$.
(For a general description of the theory of separation of variables
see Miller (1977) \cite{Miller}.)
The Euclidean distance between two points $\bfx,\bfxp\in\R^d$, expressed in a 
rotationally-invariant coordinate system, is given by
\[
\displaystyle \|\bfx-\bfxp\|=\sqrt{2RR^\prime}
\left[\chi-\cos(\phi-\phi^\prime)\right]^{1/2},
\]
where the toroidal parameter $\chi>1$, is given by \begin{equation}
\chi:=\frac{R+{R^\prime}^2
+{\displaystyle \sum_{i=3}^d(x_i-x_i^\prime)^2}
}
{\displaystyle 2RR^\prime},
\end{equation}
where $R,R^\prime\in[0,\infty)$ are defined in (\ref{rotatioanallyinvariant})
for $\bfx,\bfxp\in\R^d$.
The hypersurfaces given by $\chi>1$ equals constant are independent of
coordinate system and represent hyper-tori of revolution.

From Theorem \ref{greenpoly} we see that, apart from multiplicative constants,
the algebraic expression 
$\li_k^d:(\R^d\times\R^d)\setminus\{(\bfx,\bfx):\bfx\in\R^d\}\to\R$
of an unnormalized fundamental solution for the polyharmonic equation 
in Euclidean space $\R^d$ for $d$ even, $k\ge d/2,$ is given by
\begin{equation}
\li_k^d(\bfx,\bfxp):=\|\bfx-\bfxp\|^{2k-d}\left(\log\|\bfx-\bfxp\|-\beta_{k-d/2,d}\right).
\label{liunderscorekd}
\end{equation}
By expressing $\li_k^d$ in a rotationally-invariant coordinate system 
(\ref{rotatioanallyinvariant}) we obtain
\begin{eqnarray}
\li_k^d(\bfx,\bfxp)&=&\left(2RR^\prime\right)^p\left[\frac12\log
\left(2RR^\prime\right)-\beta_{p,d}\right]
\left[\chi-\cos(\phi-\phi^\prime) \right]^p\nonumber\\[2pt]
&&+\frac12\left(2RR^\prime\right)^p
\left[\chi-\cos(\phi-\phi^\prime) \right]^p
\log\left[\chi-\cos(\phi-\phi^\prime) \right],
\label{logfourseries}
\end{eqnarray}
where $p=k-d/2\in\N_0$.
For the polyharmonic equation in even-dimensional Euclidean space $\R^d$ with 
$1\le k \le d/2-1,$ apart from multiplicative constants, the algebraic 
expression for an unnormalized fundamental solution of the polyharmonic
equation $\hii_k^d:(\R^d\times\R^d)\setminus\{(\bfx,\bfx):\bfx\in\R^d\}\to\R$ 
is given by
\[
\hii_k^d(\bfx,\bfxp):=\|\bfx-\bfxp\|^{2k-d}.
\]
By expressing $\hii_k^d$ in a rotationally-invariant coordinate system we obtain
\begin{equation}
\hii_k^d(\bfx,\bfxp)=\left(2RR^\prime\right)^{-q}
\left[\chi-\cos(\phi-\phi^\prime) \right]^{-q},
\label{unlogfourseries}
\end{equation}
where $q=2k-d$.

By examining (\ref{logfourseries}) and (\ref{unlogfourseries}), we see that for 
computation of Fourier expansions
about the azimuthal separation angle $(\phi-\phi^\prime)$ of $\li_k^d$ and $\hii_k^d$, all that is required is to compute 
the Fourier cosine series for the following three functions
$f_\chi,h_\chi:\R\to(0,\infty)$ and $g_\chi:\R\to\R$ defined as
\begin{eqnarray*}
&f_\chi(\psi):=\left(\chi-\cos\psi\right)^p,&\\[0.2cm]
&g_\chi(\psi):=\left(\chi-\cos\psi\right)^p\log\left(\chi-\cos\psi\right),\qquad\mbox{and}&\\[0.2cm]
&h_\chi(\psi):=\left(\chi-\cos\psi\right)^{-q},&
\end{eqnarray*}
where $p\in\N_0$, $q\in\N$ and $\chi>1$ is a fixed parameter.

The Fourier series of $f_\chi$ 
is given in Cohl \& Dominici (2010) \cite{CohlDominici}
(cf.~(4.4) therein), namely
\begin{equation}
(z-\cos\psi)^p=
(z^2-1)^{p/2}
\sum_{n=0}^p
\epsilon_n\cos(n\psi)
\frac{(-p)_n(p-n)!}{(p+n)!}
P_{p}^{n}\left(\frac{z}{\sqrt{z^2-1}}\right),
\label{integer}
\end{equation}
where the Neumann factor $\epsilon_n=2-\delta_{n,0}$ commonly occurs
in Fourier series, $\delta_{n,0}$ is the Kronecker delta, and
\[
(z)_n:=\prod_{i=1}^n(z+i-1),
\]
for $z\in\C$ and $n\in\N_0$, is the Pochhammer symbol (rising factorial).
We have used Whipple's formula in (\ref{integer})
(see for instance, (8.2.7) in Abramowitz \& Stegun (1972) \cite{Abra}) 
to convert the associated Legendre function of the second kind $Q_\nu^\mu:(1,\infty)\to\C$ 
appearing in \cite{CohlDominici} to 
the associated Legendre function of the first kind $P_\nu^\mu:(1,\infty)\to\R$.
The associated Legendre function of the first kind can be defined using the 
Gauss hypergeometric function, namely (Magnus, Oberhettinger \& Soni (1966) \cite{MOS}, p.~153)
\[
P_\nu^\mu(z):=\frac{1}{\Gamma(1-\mu)}\left(\frac{z+1}{z-1}\right)^{\mu/2}
{}_2F_1\left(-\nu,\nu+1;1-\mu;\frac{1-z}{2}\right).
\]
The Gauss hypergeometric function ${}_2F_1:\C\times\C\times(\C\setminus-\N_0)\times
\{z\in\C:|z|<1\}\to\C$ can be defined 
in terms of the following infinite series
\[
{}_2F_1(a,b;c;z)=\sum_{n=0}^\infty \frac{(a)_n(b)_n}{n!(c)_n}z^n
\]
(see for instance Chapter 15 in Olver {\it et al.}~(2010) \cite{NIST}).

The Fourier series of $h_\chi$ is given in
Cohl \& Dominici (2010) \cite{CohlDominici} 
(Whipple formula 
(8.2.7) in Abramowitz \& Stegun (1972) \cite{Abra}
and cf.~(4.5) therein), namely
\begin{equation}
\hspace{-0.0cm}\frac{1}{(z-\cos\psi)^q}=
\frac{(z^2-1)^{-q/2}}{(q-1)!}
\sum_{n=0}^\infty\epsilon_n\cos(n\psi)(n+q-1)! P_{q-1}^{-n}
\left(\frac{z}{\sqrt{z^2-1}}\right),
\label{azimuthalfourierseriesofoneoverq}
\end{equation}
where $q\in\N$.
Since the Fourier series of $h_\chi$ is computed in Cohl \& 
Dominici (2010) \cite{CohlDominici}, we understand how to compute 
Fourier expansions of $\hii_k^d$ (\ref{unlogfourseries})
in separable rotationally-invariant coordinate systems.
In order to compute Fourier expansion of $\li_k^d$ (\ref{logfourseries})
in separable rotationally-invariant coordinate systems, all that remains is to determine 
the Fourier series of $g_\chi$. This is the goal of the next two sections.

\section{Algebraic approach to the logarithmic Fourier series}
\label{AlgebraicapproachtologarithmicFourierseries}

Since $\chi>1$, one may make the substitution $\chi=\cosh\eta$ to evaluate 
the Fourier series of $g_\chi$.
For instance, it is given in the form of 
$(\cosh\eta-\cos\psi)^p\log(\cosh\eta-\cos\psi),$ where $p\in\N_0$. 
For $p=0$ the result is well-known
(see for instance Magnus, Oberhettinger \& Soni (1966) \cite{MOS}, p.~259)
\begin{equation}
\log(\cosh\eta-\cos\psi)=\eta-\log{2}-2\sum_{n=1}^\infty \frac{e^{-n\eta}}{n} \cos(n\psi),
\label{logcoshcos}
\end{equation}
which as we will see, should be compared with
(\ref{azimuthalfourierseriesofoneoverq}) for $q=1$, namely
\begin{equation}
\displaystyle \frac{1}{\cosh\eta-\cos\psi}=
\frac{1}{\sinh\eta}
\sum_{n=0}^\infty
\epsilon_n
\cos(n\psi)
e^{-n\eta}.
\label{genint111}
\end{equation}
Note that for $\eta>0$ we may write $e^\eta$ and therefore $\eta$ as a function of
$\cosh\eta$ since
$\sinh\eta=\sqrt{\cosh^2\eta-1},$
$e^\eta=\cosh\eta+\sqrt{\cosh^2\eta-1},$
and therefore
$\eta=\log\left(\cosh\eta+\sqrt{\cosh^2\eta-1}\right).$
\noindent Now examine the $p=1$ case for $g_\chi$.  If we multiply both sides of 
(\ref{logcoshcos}) by $(\cosh\eta-\cos\psi)$ and take advantage of the formula 
\begin{equation}
\cos(n\psi)\cos\psi=\frac12\Bigl\{\cos[(n+1)\psi]+\cos[(n-1)\psi]\Bigr\},
\label{prod2cos}
\end{equation}
then we have 
\begin{eqnarray}
\hspace{-1.5cm}(\cosh\eta-\cos\psi)\log(\cosh\eta-\cos\psi)&=&
(\eta-\log 2)\cosh\eta\nonumber\\[0.2cm]
&&{}\hspace{-6cm}-(\eta-\log 2)\cos\psi
-2\cosh\eta\sum_{n=1}^\infty \frac{e^{-n\eta}}{n}\cos(n\psi)\nonumber\\[0.2cm]
&&{}\hspace{-6cm}+\sum_{n=1}^\infty \frac{e^{-n\eta}}{n}\cos[(n+1)\psi]
+\sum_{n=1}^\infty \frac{e^{-n\eta}}{n}\cos[(n-1)\psi].
\end{eqnarray}
Collecting the contributions to the Fourier cosine series, we obtain
\begin{eqnarray}
\hspace{-1.0cm}(\cosh\eta-\cos\psi)\log(\cosh\eta-\cos\psi)&=&(1+\eta-\log 2)\cosh\eta
\nonumber\\[0.1cm]
&&{}\hspace{-2.2cm}-\sinh\eta+\cos\psi\left(\log 2-1-\eta-\frac12 e^{-2\eta}\right)\nonumber\\[0.2cm]
&&{}\hspace{-2.2cm}+2\sum_{n=2}^\infty \frac{e^{-n\eta}\cos n\psi}{n(n^2-1)}(\cosh\eta+n\sinh\eta).
\label{logfourierp1}
\end{eqnarray}
If we compare (\ref{logfourierp1}) with (\ref{azimuthalfourierseriesofoneoverq}) for $q=1$, namely
\begin{equation}
\displaystyle \frac{1}{(\cosh\eta-\cos\psi)^2}=
\frac{1}{\sinh^3\eta}
\sum_{n=0}^\infty
\epsilon_n
\cos(n\psi)
e^{-n\eta}
(\cosh\eta+n\sinh\eta),
\label{genint2}
\end{equation}
we notice that the factor $(\cosh\eta+n\sinh\eta)$ appears in both series.

For $p=2$ in $g_\chi$, we use (\ref{prod2cos}) and similarly have
\begin{eqnarray}
\hspace{-0.0cm}(\cosh\eta-\cos\psi)^2\log(\cosh\eta-\cos\psi)&&
\nonumber\\[0.2cm]
&&{}\hspace{-5.5cm}=(\eta-\log 2)\cosh^2\eta
-2(\eta-\log 2)\cosh\eta\cos\psi\nonumber\\[0.2cm]
&&{}\hspace{-5.3cm}
+(\eta-\log 2)\cos^2\psi
-(2\cosh^2\eta+1)\sum_{n=1}^\infty \frac{e^{-n\eta}}{n}
\cos(n\psi)\nonumber
\nonumber\\[0.2cm]
&&{}\hspace{-5.3cm}
+2\cosh\eta\sum_{n=1}^\infty \frac{e^{-n\eta}}{n}\cos[(n+1)\psi]
+2\cosh\eta\sum_{n=1}^\infty \frac{e^{-n\eta}}{n}\cos[(n-1)\psi]
{}\nonumber\\[0.2cm]
&&{}\hspace{-5.3cm}
-\frac12\sum_{n=1}^\infty \frac{e^{-n\eta}}{n}\cos[(n+2)\psi]
-\frac12\sum_{n=1}^\infty \frac{e^{-n\eta}}{n}\cos[(n-2)\psi].\nonumber
\end{eqnarray}
If we collect the contributions of the Fourier cosine series, we obtain
\begin{eqnarray}
\hspace{-1cm}(\cosh\eta-\cos\psi)^2\log(\cosh\eta-\cos\psi)&=&
(\eta-\log 2)\left(\cosh^2\eta+\frac12\right)
\nonumber\\[0.1cm]
&&{}\hspace{-6.2cm}
+2\cosh\eta\,e^{-\eta}-\frac14 e^{-2\eta}
+\Biggl[-2(\eta-\log 2)\cosh\eta
-\left(2\cosh^2\eta+\frac32\right)e^{-\eta}\nonumber\\[0.1cm]
&&{}\hspace{-6.2cm} +\cosh\eta\,e^{-2\eta}
-\frac16 e^{-3\eta}
\Biggr]\cos\psi
+\Biggl[\frac12(\eta-\log 2) +2\cosh\eta\,e^{-\eta}
\nonumber\\[0.1cm]
&&{}\hspace{-6.2cm}
-\frac12 (2\cosh^2\eta+1)e^{-2\eta}
+\frac23 \cosh\eta\,e^{-3\eta}-\frac18e^{-4\eta}\Biggr]\cos 2\psi\nonumber\\[0.2cm]
&&{}\hspace{-6.2cm}-4\sum_{n=3}^\infty \frac{e^{-n\eta}\cos(n\psi)}{n(n^2-1)(n^2-4)}
\left[(n^2-1)\sinh^2\eta+3n\sinh\eta\cosh\eta+3\cosh^2\eta\right].
\label{logfourierp2}
\end{eqnarray}
By comparing (\ref{logfourierp2}) with (\ref{azimuthalfourierseriesofoneoverq}) for $q=2$, namely
\begin{eqnarray}
\hspace{0.0cm}\displaystyle\hspace{-1truecm}\frac{1}{(\cosh\eta-\cos\psi)^3}&=&
\frac{1}{2\sinh^5\eta}
\nonumber\\[0.2cm]
&&{}\hspace{-2.7cm}
\times\sum_{n=0}^\infty
\epsilon_n
\cos(n\psi)
e^{-n\eta}
\left[(n^2-1)\sinh^2\eta
+3n\sinh\eta\cosh\eta+3\cosh^2\eta\right],
\label{genint3}
\end{eqnarray}
then we notice that the factor
$\left((n^2-1)\sinh^2\eta
+3n\sinh\eta\cosh\eta+3\cosh^2\eta\right)$ appears in both series.  We will 
demonstrate in \S\ref{Comparison}, 
why the identification mentioned in 
(\ref{genint2}) and (\ref{genint3}) occurs.

\medskip
This algebraic approach for determining the Fourier series of $g_\chi$ will
now be generalized.
By starting with (\ref{logcoshcos})
and repeatedly multiplying by factors of $(\cosh\eta-\cos\psi)$, we see that the general 
Fourier series of $g_\chi$ can be given in terms of a sequence of polynomials 
$R_p^k:(1,\infty)\to\R$, with $p\in\N_0$ and $k\in\Z$, as
\begin{eqnarray}
\hspace{-1cm}(\cosh\eta-\cos\psi)^p\log(\cosh\eta-\cos\psi)&=&
(\eta-\log 2)(\cosh\eta-\cos\psi)^p\nonumber\\[0.1cm]
&&\hspace{-4.5cm}{}+2\sum_{k=-p}^{p}(-1)^{k+1}R_p^k(\cosh\eta)
\sum_{n=1}^\infty
\frac{e^{-n\eta}}{n}
\cos[(n+k)\psi].
\label{coshplogcoshp1}
\end{eqnarray}

\noindent We will refer to $R_p^k$ as logarithmic polynomials
with argument $x=\cosh\eta$ 
(in our notation $p$ and $k$ are both indices) (See 
Appendix \ref{Thelogarithmicpolynomials} for a description of some of the
properties of the logarithmic polynomials).

The double sum in (\ref{coshplogcoshp1}) is simplified by making the replacement
$n+k\mapsto n$.  It then follows that the resulting double sum naturally breaks 
into two disjoint regions, one triangular
\[
{\mathcal A}:=\{(k,n):-p\le k \le p-1,\ k+1\le n\le p\},
\]
with $p(2p+1)$ terms and the other infinite rectangular
\[
{\mathcal B}:=\{(k,n):-p\le k \le p,\ p+1<n<\infty\}.
\]
By rearranging the order of the $k$ and $n$ summations
in (\ref{coshplogcoshp1}), we derive
\begin{eqnarray}
\hspace{-0.4cm}\displaystyle (\cosh\eta-\cos\psi)^p\log(\cosh\eta-\cos\psi)
&=&(\eta-\log 2)(\cosh\eta-\cos\psi)^p\nonumber\\[0.2cm]
&&\hspace{-5.8cm}\displaystyle 
+\sum_{n=0}^{p}\cos(n\psi)e^{-n\eta}
\,\mathfrak{r}_{n,p}^{-p,n-1}(\cosh\eta)
+\sum_{n=1}^{p-1}\cos(n\psi)e^{n\eta}
\,\mathfrak{r}_{-n,p}^{-p,-n-1}(\cosh\eta)
\nonumber\\[0.0cm]
&&\hspace{-5.8cm}\displaystyle
+\sum_{n=p+1}^\infty\frac{\cos(n\psi)e^{-n\eta}}
{n(n^2-1)\cdots(n^2-p^2)}
\,\Re_{n,p}(\cosh\eta),
\label{coshplogcoshp2}
\end{eqnarray}
\noindent where $
\mathfrak{r}_{n,p}^{k_1,k_2},
\Re_{n,p}
:(1,\infty)\to\R$ 
are defined as
\[
\mathfrak{r}_{n,p}^{k_1,k_2}(\cosh\eta):=
2\sum_{k=k_1}^{k_2}
\frac{
(-1)^{k+1}e^{k\eta}R_p^k(\cosh\eta)
}{n-k},
\]
and
\[
\Re_{n,p}(\cosh\eta):=\frac{(n+p)!}{(n-p-1)!}
\mathfrak{r}_{n,p}^{-p,p}(\cosh\eta),\quad (n\ge p+1),
\]
respectively.  
We can also write the Fourier series directly in terms of the logarithmic polynomials
$R_p^k$ as follows
\begin{eqnarray}
\hspace{-0.5cm}(\cosh\eta-\cos\psi)^p\log(\cosh\eta-\cos\psi)&=&(\eta-\log 2)(\cosh\eta-\cos\psi)^p\nonumber\\[0.2cm]
&&\hspace{-5.3cm}+2\sum_{n=0}^p\cos(n\psi)e^{-n\eta}\sum_{k=-p}^{n-1}
\frac{(-1)^{k+1}e^{k\eta}R_p^k(\cosh\eta)}{n-k}\nonumber\\[0.2cm]
&&\hspace{-5.3cm}-2\sum_{n=1}^{p-1}\cos(n\psi)e^{n\eta}\sum_{k=-p}^{-n-1}
\frac{(-1)^{k+1}e^{k\eta}R_p^k(\cosh\eta)}{n+k}\nonumber\\[0.2cm]
&&\hspace{-5.3cm}+2\sum_{n=p+1}^\infty\cos(n\psi)e^{-n\eta}\sum_{k=-p}^{p}
\frac{(-1)^{k+1}e^{k\eta}R_p^k(\cosh\eta)}{n-k}.
\label{cosetalogintermsofRpk}
\end{eqnarray}
Note that by using (\ref{integer}), then we can express $(\cosh\eta-\cos\psi)^p$ as
a Fourier series, namely
\begin{equation}
(\cosh\eta-\cos\psi)^p=\sinh^p\eta
\sum_{n=0}^p
\epsilon_n\cos(n\psi)
\frac{(-p)_n(p-n)!}{(p+n)!}
P_{p}^{n}(\coth\eta),
\label{coshetacospsip}
\end{equation}
where $p\in\N_0$.

If we define the function 
$\mathfrak{P}_{n,p}:(1,\infty)\to\R$ such that
\begin{eqnarray}
\mathfrak{P}_{n,p}(\cosh\eta):=
\left\{ \begin{array}{ll}
\displaystyle \mathfrak{r}_{0,p}^{-p,-1}(\cosh\eta) 
\hspace{4.578cm}\mathrm{if}\ n=0, \nonumber \\[0.2cm] 
\displaystyle D_{n,p}(\cosh\eta) 
+E_{n,p}(\cosh\eta) 
\hspace{2.19cm}\mathrm{if}\ 1\le n\le p-1,\nonumber \\[0.2cm]
\displaystyle e^{-p\eta}\mathfrak{r}_{p,p}^{-p,p-1}(\cosh\eta)
\hspace{3.68cm}\mathrm{if}\ 1\le n=p,\nonumber\\[0.2cm]
\displaystyle\frac{e^{-n\eta}}
{(n^2-p^2)\cdots(n^2-1)n}
\,\Re_{n,p}(\cosh\eta)\hspace{0.75cm}\mathrm{if}\ n\ge p+1,
\end{array} \right.
\end{eqnarray}
where $D_{n,p},E_{n,p}:(1,\infty)\to\R$ are defined as
\[
D_{n,p}(\cosh\eta)=
\left\{ \begin{array}{ll}
\displaystyle e^{n\eta}\mathfrak{r}_{-n,p}^{-p,-n-1}(\cosh\eta) &\qquad\mathrm{if}\ p\ge 2, \nonumber \\[0.2cm]
\displaystyle 0 &\qquad\mathrm{if}\ p\in\{0,1\}, \nonumber
\end{array} \right.
\]
and 
\[
E_{n,p}(\cosh\eta)=
\left\{ \begin{array}{ll}
\displaystyle e^{-n\eta}\mathfrak{r}_{n,p}^{-p,n-1}(\cosh\eta) &\qquad\mathrm{if}\ p\ge 1, \nonumber \\[0.2cm]
\displaystyle 0 &\qquad\mathrm{if}\ p=0, \nonumber
\end{array} \right.
\]
respectively, then we can write (\ref{coshplogcoshp2}) as 
\begin{eqnarray}
\displaystyle\hspace{-1cm}(\cosh\eta-\cos\psi)^p\log(\cosh\eta-\cos\psi)&=&
(\eta-\log 2)(\cosh\eta-\cos\psi)^p\nonumber\\[0.2cm]
&&\hspace{0cm}+\sum_{n=0}^\infty\cos(n\psi)
\mathfrak{P}_{n,p}(\cosh\eta).\nonumber
\end{eqnarray}
In fact, if we use (\ref{coshetacospsip}), then we can express $g_\chi$ as
\begin{eqnarray}
\hspace{0.2cm}\displaystyle(\cosh\eta-\cos\psi)^p\log(\cosh\eta-\cos\psi)&=&
\sum_{n=0}^\infty\cos(n\psi)
\mathfrak{P}_{n,p}(\cosh\eta)\nonumber\\[0.2cm]
&&{}\hspace{-6.2cm}
+(\eta-\log 2)\sinh^p\eta
\sum_{n=0}^p
\epsilon_n\cos(n\psi)
\frac{(-p)_n(p-n)!}{(p+n)!}
P_{p}^{n}(\coth\eta).\nonumber
\end{eqnarray}
Furthermore, if we define $\mathfrak{Q}_{n,p}:(1,\infty)\to\R$ as
\[
\mathfrak{Q}_{n,p}(\cosh\eta):=\mathfrak{P}_{n,p}(\cosh\eta)
+\frac{\epsilon_n(-p)_n(p-n)!}{(p+n)!}
(\eta-\log 2)\sinh^p\eta
P_{p}^{n}(\coth\eta),
\]
then we can express $g_\chi$ as
\begin{equation}
\displaystyle(\cosh\eta-\cos\psi)^p\log(\cosh\eta-\cos\psi)=
\sum_{n=0}^\infty\cos(n\psi)
\mathfrak{Q}_{n,p}(\cosh\eta).
\label{Qnpformoflogexpansion}
\end{equation} 

\section{Limit derivative approach to the logarithmic Fourier series}
\label{LimitderivativeapproachtologarithmicFourierseries}

We now use a second approach to compute the Fourier series for a logarithmic 
fundamental solution of the polyharmonic equation
(\ref{liunderscorekd}).
We would like to match our results to 
the computations in \S \ref{AlgebraicapproachtologarithmicFourierseries}, which 
clearly demonstrate different behaviors for the two regimes,
$0\le n\le p$ and $n\ge p+1$.  
By applying the identity
\begin{equation}
(\cosh\eta-\cos\psi)^p\log(\cosh\eta-\cos\psi)=
\lim_{\nu\rightarrow 0} \frac{\partial}{\partial\nu}
(\cosh\eta-\cos\psi)^{\nu+p},
\label{logconsequence}
\end{equation}
where $p\in\N_0$, to
\begin{eqnarray}
\hspace{-0.4cm}(\cosh\eta-\cos\psi)^\nu&=&
\sinh^\nu\eta\sum_{n=0}^\infty
\frac{(-1)^n\epsilon_n\cos(n\psi)}{(\nu+1)_n}
P_{\nu}^{n}(\coth\eta),
\label{expansiondeq2}
\end{eqnarray}
where 
$\nu\in\C\setminus-\N$ (cf.~(3.11b) in Cohl \& Dominici (2010) \cite{CohlDominici}),
one can compute the Fourier cosine series of $(\cosh\eta-\cos\psi)^p\log(\cosh\eta-\cos\psi)$,
provided availability of the necessary parameter derivatives.

Applying (\ref{logconsequence}) to (\ref{expansiondeq2}), we obtain
\begin{eqnarray*}
(\cosh\eta-\cos\psi)^p\log(\cosh\eta-\cos\psi)&&\\[0.2cm]
&&\hspace{-5.0cm}=
\left[\lim_{\nu\to 0}\frac{\partial}{\partial\nu}\sinh^{\nu+p}\eta\right]
\sum_{n=0}^\infty\frac{(-1)^n\epsilon_n\cos(n\psi)}{(p+1)_n}P_p^n(\coth\eta)\\[0.2cm]
&&\hspace{-5.0cm}+\sinh^p\eta
\sum_{n=0}^\infty (-1)^n\epsilon_n\cos(n\psi)
\left[\lim_{\nu\to 0}\frac{\partial}{\partial\nu}\frac{1}{(\nu+p+1)_n}\right]
P_p^n(\coth\eta)\\[0.2cm]
&&\hspace{-5.0cm}+\sinh^p\eta
\sum_{n=0}^\infty \frac{(-1)^n\epsilon_n\cos(n\psi)}{(p+1)_n}
\left[\lim_{\nu\to 0}\frac{\partial}{\partial\nu}
P_{\nu+p}^n(\coth\eta)
\right].
\end{eqnarray*}
Note that for $p\in\Z$ and $n\in\N_0$, the associated Legendre function of 
the first kind $P_p^n$ vanishes if $n\ge p+1$.
This is easily seen using the Rodrigues-type formula (cf.~(14.7.11) in
Olver {\it et al.} (2010) \cite{NIST})
\[
P_p^n(z)=(z^2-1)^{n/2}\frac{d^nP_p(z)}{dz^n},
\]
and the fact that $P_p(z)$ is a polynomial in $z$ of degree $p$.
The derivatives are given as follows:
\begin{equation}
\lim_{\nu\to 0}\frac{\partial}{\partial\nu}\sinh^{\nu+p}\eta
=\sinh^p\eta\log\sinh\eta,
\label{der1}
\end{equation}
\begin{equation}
\lim_{\nu\to 0}\frac{\partial}{\partial\nu}\frac{1}{(\nu+p+1)_n}
=\frac{p!}{(p+n)!}\left[\psi(p+1)-\psi(p+1+n)\right],
\label{der2}
\end{equation}
and
\begin{equation}
\lim_{\nu\to 0}\frac{\partial}{\partial\nu}
P_{\nu+p}^n(\coth\eta)
=\left[\frac{\partial}{\partial\nu}P_\nu^n(\coth\eta)\right]_{\nu=p},
\label{ddlfotfk}
\end{equation}
where $\psi:\C\setminus-\N_0\to\C$ is the digamma function
defined in terms of the 
derivative of the gamma function
\[
\frac{d}{dz}\Gamma(z)=:\psi(z)\Gamma(z)
\]
(see for instance (5.2.2) in Olver {\it et al.} (2010) \cite{NIST}).
The degree-derivative of the associated Legendre function of the first kind
in (\ref{ddlfotfk}) is determined 
using (\ref{limitderivdegreepnum}) and (\ref{limitderivdegreepnum2}).
By collecting terms and using (\ref{der1}), (\ref{der2}), and (\ref{ddlfotfk}), we obtain
\begin{eqnarray}
\hspace{-0.5cm}(\cosh\eta-\cos\psi)^p\log(\cosh\eta-\cos\psi)&=&
(\eta-\log 2)(\cosh\eta-\cos\psi)^p\nonumber\\[0.2cm]
&&\hspace{-6.0cm}+p!\sinh^p\eta\sum_{n=0}^p
\frac{(-1)^n\epsilon_n\cos(n\psi)}{(p+n)!}\nonumber\\[0.2cm]
&&\hspace{-5.0cm}\times
\left[
2\psi(2p+1)-\psi(p+1+n)-\psi(p+1-n)
\right]P_p^n(\coth\eta)\nonumber\\[0.2cm]
&&\hspace{-6.0cm}
+(-1)^pp!\sinh^p\eta\sum_{n=0}^{p-1}\frac{\epsilon_n\cos(n\psi)}{(p+n)!}\nonumber\\[0.2cm]
&&\hspace{-5.0cm}
\times\sum_{k=0}^{p-n-1}\frac
{(-1)^k(2n+2k+1)\left[1+\frac{k!(p+n)!}{(2n+k)!(p-n)!}\right]}
{(p-n-k)(p+n+k+1)}
P_{n+k}^n(\coth\eta)\nonumber\\[0.2cm]
&&\hspace{-6.0cm}
+2(-1)^pp!\sinh^p\eta\sum_{n=1}^{p}\frac{(-1)^n\cos(n\psi)}{(p-n)!}
\nonumber\\[0.2cm] &&\hspace{-5.0cm}\times
\sum_{k=0}^{n-1}\frac
{(-1)^k(2k+1)}
{(p-k)(p+k+1)}
P_{k}^{-n}(\coth\eta)\nonumber\\[0.2cm]
&&\hspace{-6.0cm}
+2(-1)^{p+1}p!\sinh^p\eta\sum_{n=p+1}^\infty\cos(n\psi)(n-p-1)!P_p^{-n}(\coth\eta).
\label{limitderivativeapproach}
\end{eqnarray}

\section{Comparison of the two approaches}
\label{Comparison}

The limit derivative approach presented in
\S\ref{LimitderivativeapproachtologarithmicFourierseries}
might be considered, of the two methods, preferred for 
computing the azimuthal Fourier series for a logarithmic fundamental solution 
of the polyharmonic equation.  This is because it produces azimuthal Fourier
coefficients in terms of the well-known special functions, associated Legendre 
functions.  On the other hand, the algebraic approach presented in 
\S\ref{AlgebraicapproachtologarithmicFourierseries} produces results in 
terms of the
two-parameter 
logarithmic polynomials $R_p^k(x)$. As far as the author is aware, these
polynomials are previously unencountered in the literature.  By comparison 
of the two
approaches we see how the logarithmic polynomials $R_p^k(x)$ (potentially 
a new type of {\it special function}) are intimately related to the associated 
Legendre functions.  In this section we make this comparison concrete.  We
should also mention that the following comparison equations resolve to become 
quite complicated as $p$ increases, and they have been checked for $0\le p\le 10$
using Mathematica with the assistance of an algorithm generated using 
(\ref{algo1}) and (\ref{algo2}) from Appendix \ref{Thelogarithmicpolynomials}.
\smallskip

By equating the Fourier coefficients using the two approaches we can obtain
summation formulae which are satisfied by the logarithmic polynomials.  
For $n=0$ and $p\ge 1$ we have 
\begin{eqnarray*}
\hspace{-0.2cm}\sum_{k=1}^{p}\frac{(-1)^{k+1}e^{-k\eta}R_p^k(\cosh\eta)}{k}&&\\[0.2cm]
&&{}\hspace{-3.5cm}
=\sinh^p\eta\left[\psi(2p+1)-\psi(p+1)\right]P_p(\coth\eta)\\[0.2cm]
&&{}\hspace{-2.5cm}
+(-1)^p\sinh^p\eta\sum_{k=0}^{p-1}\frac{(-1)^k(2k+1)}{(p-k)(p+k+1)}P_k(\coth\eta),
\end{eqnarray*}
for $1\le n\le p-1$, $p\ge 2$ we derive
\begin{eqnarray*}
\hspace{-0.4cm}\sum_{k=-p}^{n-1}\frac{(-1)^{k+1}e^{k\eta}R_p^k(\cosh\eta)}{n-k}
+e^{2n\eta}\sum_{k=-p}^{-n-1}\frac{(-1)^{k}e^{k\eta}R_p^k(\cosh\eta)}{n+k}
&=&p!e^{n\eta}\sinh^p\eta\\[0.2cm]
&&{}\hspace{-10.7cm}
\times\Biggl\{
\frac{(-1)^n}{(p+n)!}\left[2\psi(2p+1)-\psi(p+1+n)-\psi(p+1-n)\right]P_p^n(\coth\eta)\\[0.2cm]
&&{}\hspace{-9.9cm}
+\frac{(-1)^p}{(p+n)!}\sum_{k=0}^{p-n-1}\frac{(-1)^k(2n+2k+1)
\left[1+\frac{k!(p+n)!}{(2n+k)!(p-n)!}\right]
}{(p-n-k)(p+n+k+1)}
P_{n+k}^n(\coth\eta)\\[0.2cm]
&&{}\hspace{-9.9cm}
+\frac{(-1)^{p+n}}{(p-n)!}\sum_{k=0}^{n-1}\frac{(-1)^k(2k+1)}{(p-k)(p+k+1)}P_k^{-n}(\coth\eta)
\Biggr\},
\end{eqnarray*}
for $n=p$, $p\ge 1$ we obtain 
\begin{eqnarray*}
\hspace{-0.2cm}\sum_{k=-p}^{p-1}\frac{(-1)^{k+1}e^{k\eta}R_p^k(\cosh\eta)}{p-k}
&=&p!e^{p\eta}\sinh^p\eta\\[0.2cm]
&&{}\hspace{-3.2cm}
\times\Biggl\{
\frac{(-1)^p}{(2p)!}\left[\psi(2p+1)-\psi(1)\right]P_p^p(\coth\eta)\\[0.2cm]
&&{}\hspace{-0.7cm}
+\sum_{k=0}^{p-1}\frac{(-1)^k(2k+1)}{(p-k)(p+k+1)}P_k^{-p}(\coth\eta)
\Biggr\},
\end{eqnarray*}
and for $n\ge p+1$, $p\ge 0,$ we see that
\begin{eqnarray}
\hspace{-1.0cm}\sum_{k=-p}^{p}\frac{(-1)^{k+1}e^{k\eta}R_p^k(\cosh\eta)}{n-k}&&\nonumber\\[0.2cm]
&&\hspace{-2.5cm}=(-1)^{p+1}p!(n-p-1)!e^{n\eta}\sinh^p\eta P_p^{-n}(\coth\eta).
\label{ngepp1}
\end{eqnarray}

We now have closed-form expressions for the
finite terms given by
(\ref{coshplogcoshp2}),
(\ref{Qnpformoflogexpansion}) and
(\ref{coshplogcoshp2}) in terms of associated 
Legendre functions of the first kind.
We also have a proof 
of the correspondence for the ``ending'' function 
$\Re_{n,q-1}$ 
mentioned in \S\ref{AlgebraicapproachtologarithmicFourierseries}.
Through (\ref{ngepp1}), the function $\Re_{n,p}$ 
(cf.~(\ref{logfourierp1}), (\ref{genint2}),
(\ref{logfourierp2}),
(\ref{genint3}), and
(\ref{coshplogcoshp2})) 
is directly related to
the associated Legendre function of the first kind, namely
\[
\Re_{n,p}(\cosh\eta)=
2(-1)^{p+1}p!(p+n)!e^{n\eta}\sinh^p\eta P_p^{-n}(\coth\eta).
\]
Therefore through (\ref{azimuthalfourierseriesofoneoverq}) we have
\[
\frac{1}{(\cosh\eta-\cos\psi)^q}=
\frac{(-1)^q}{2[(q-1)!]^2\sinh^{2q-1}\eta}
\sum_{n=0}^\infty
\epsilon_n
\cos(n\psi)
e^{-n\eta}
\Re_{n,q-1}(\cosh\eta),
\]
which demonstrates the correspondences which was mentioned near
(\ref{genint111}), (\ref{genint2}), and (\ref{genint3}) for (\ref{coshplogcoshp2}) and
(\ref{limitderivativeapproach}).

\section{Fourier expansion for a logarithmic fundamental solution of the polyharmonic equation}
\label{FourierExpansionforlogfundsolofpolyeq}

Now that we have computed the Fourier series for $g_\chi$, namely 
(\ref{coshplogcoshp2})
(cf.~(\ref{Qnpformoflogexpansion}))
and (\ref{limitderivativeapproach}), we are in a position 
to compute the azimuthal Fourier series for a logarithmic fundamental
solution of the polyharmonic equation (\ref{liunderscorekd}).



For instance, using (\ref{logfourseries}), (\ref{coshetacospsip}) 
and (\ref{limitderivativeapproach}), we have
\begin{eqnarray*}
\hspace{0.1cm}\li_k^d(\bfx,\bfxp)&=&
\frac12\left(2RR^\prime\right)^p\left[\log\left(RR^\prime\right)
+\eta-\beta_{p,d}\right](\chi^2-1)^{p/2}\\[0.2cm]
&&\hspace{-0.0cm}\times \sum_{n=0}^p\epsilon_n\cos[n(\phi-\phi^\prime)]
\frac{(-p)_n(p-n)!}{(p+n)!}P_p^n\Biggl(\frac{\chi}{\sqrt{\chi^2-1}}\Biggr)\\[0.2cm]
&&\hspace{-1.8cm}+\frac12(2RR^\prime)^pp!(\chi^2-1)^{p/2}
\sum_{n=0}^p
\frac{(-1)^n\epsilon_n\cos[n(\phi-\phi^\prime)]}{(p+n)!}\nonumber\\*[0.2cm]
&&\hspace{-0.0cm}\times
\left[
2\psi(2p+1)-\psi(p+1+n)-\psi(p+1-n)
\right]P_p^n\Biggl(\frac{\chi}{\sqrt{\chi^2-1}}\Biggr)\nonumber\\[0.2cm]
&&\hspace{-1.8cm}
+\frac12(2RR^\prime)^pp!(\chi^2-1)^{p/2}
\sum_{n=0}^{p-1}\frac{\epsilon_n\cos[n(\phi-\phi^\prime)]}{(p+n)!}\nonumber\\[0.2cm]
&&\hspace{-0.0cm}
\times\sum_{k=0}^{p-n-1}\frac
{(-1)^k(2n+2k+1)\left[1+\frac{k!(p+n)!}{(2n+k)!(p-n)!}\right]}
{(p-n-k)(p+n+k+1)}
P_{n+k}^n\Biggl(\frac{\chi}{\sqrt{\chi^2-1}}\Biggr)\nonumber\\[0.2cm]
&&\hspace{-1.8cm}
+(2RR^\prime)^p(-1)^pp!(\chi^2-1)^{p/2}
\sum_{n=1}^{p}\frac{(-1)^n\cos[n(\phi-\phi^\prime)]}{(p-n)!}
\nonumber\\[0.2cm] 
&&\hspace{-0.0cm}\times
\sum_{k=0}^{n-1}\frac
{(-1)^k(2k+1)}
{(p-k)(p+k+1)}
P_{k}^{-n}\Biggl(\frac{\chi}{\sqrt{\chi^2-1}}\Biggr)\nonumber\\[0.2cm]
&&\hspace{-1.8cm}
+(2RR^\prime)^p(-1)^{p+1}p!(\chi^2-1)^{p/2}\\[0.2cm]
&&\hspace{-0.0cm}\times\sum_{n=p+1}^\infty\cos[n(\phi-\phi^\prime)](n-p-1)!
P_p^{-n}\Biggl(\frac{\chi}{\sqrt{\chi^2-1}}\Biggr).
\end{eqnarray*}
Alternatively, using (\ref{logfourseries}), 
(\ref{integer}), and
(\ref{cosetalogintermsofRpk}) we obtain
\begin{eqnarray*}
\hspace{0.5cm}\li_k^d(\bfx,\bfxp)&=&
\frac12\left(2RR^\prime\right)^p\left[\log\left(RR^\prime\right)
+\eta-\beta_{p,d}\right](\chi^2-1)^{p/2}\\[0.2cm]
&&\hspace{2.0cm}\times\sum_{n=0}^p \epsilon_n\cos[n(\phi-\phi^\prime] \frac{(-p)_n(p-n)!}{(p+n)!}
P_{p}^{n}\left(\frac{\chi}{\sqrt{\chi^2-1}}\right),\\[0.2cm]
&&\hspace{-0.3cm}
+(2RR^\prime)^p\sum_{n=0}^p\cos[n(\phi-\phi^\prime)]
\left(\chi-\sqrt{\chi^2-1}\right)^n\\[0.2cm]
&&\hspace{2.0cm}
\times\sum_{k=-p}^{n-1}
\frac{(-1)^{k+1}
\left(\chi+\sqrt{\chi^2-1}\right)^k
R_p^k(\chi)}{n-k}\nonumber\\[0.2cm]
&&\hspace{-0.3cm}
-(2RR^\prime)^p\sum_{n=1}^{p-1}\cos[n(\phi-\phi^\prime)]
\left(\chi+\sqrt{\chi^2-1}\right)^n\\[0.2cm]
&&\hspace{2.0cm}
\times\sum_{k=-p}^{-n-1}
\frac{(-1)^{k+1}
\left(\chi+\sqrt{\chi^2-1}\right)^k
R_p^k(\chi)}{n+k}\nonumber\\[0.2cm]
&&\hspace{-0.3cm}
+(2RR^\prime)^p\sum_{n=p+1}^\infty\cos[n(\phi-\phi^\prime)]
\left(\chi-\sqrt{\chi^2-1}\right)^n\\[0.2cm]
&&\hspace{2.0cm}
\times\sum_{k=-p}^{p}
\frac{(-1)^{k+1}
\left(\chi+\sqrt{\chi^2-1}\right)^k
R_p^k(\chi)}{n-k}.\nonumber
\end{eqnarray*}

Using these formulae, we can, for instance, obtain the axisymmetric component of a 
logarithmic fundamental solution of the polyharmonic equation, namely
\begin{eqnarray*}
\hspace{-0.1cm}\left.\li_k^d(\bfx,\bfxp)\right|_{n=0}&=&\frac12(2RR^\prime)^p(\chi^2-1)^{p/2}\\[0.2cm]
&&{}\hspace{-1.0cm}\times\left[\log(RR^\prime)+\eta-\beta_{p,d}+2\psi(2p+1)-2\psi(p+1)\right]
P_{p}\Biggl(\frac{\chi}{\sqrt{\chi^2-1}}\Biggr)\nonumber\\[0.2cm]
&&{}\hspace{-1.0cm}+(2RR^\prime)^p(\chi^2-1)^{p/2}\sum_{k=0}^{p-1}
\frac{(-1)^k(2k+1)}{(p-k)(p+k+1)}
P_{k}\Biggl(\frac{\chi}{\sqrt{\chi^2-1}}\Biggr),
\end{eqnarray*}
or 
\begin{eqnarray*}
\hspace{-0.4cm}\left.\li_k^d(\bfx,\bfxp)\right|_{n=0}&=&\frac12(2RR^\prime)^p(\chi^2-1)^{p/2}
\left[\log(RR^\prime)+\eta-\beta_{p,d}\right]
P_{p}\Biggl(\frac{\chi}{\sqrt{\chi^2-1}}\Biggr)\nonumber\\[0.2cm]
&&{}\hspace{1.8cm}+(2RR^\prime)^p\sum_{k=1}^{p}
\frac{(-1)^{k+1}(\chi-\sqrt{\chi^2-1})^kR_p^k(\chi)}{k}.
\end{eqnarray*}
The above expressions for the axisymmetric component of a fundamental solution of
the polyharmonic equation is one type of expression sought after in Tsai, Chen \& Hsu (2009)
\cite{TsaiChenHsu}.

\appendices


\section{Derivatives with respect to the degree of certain integer-order
associated Legendre functions of the first kind}
\label{DerivativeswithrespecttothedegreeofP}

The derivative with respect to its degree for the associated Legendre function
of the first kind evaluated at the zero degree is given in \S 4.4.3 of
Magnus, Oberhettinger \& Soni (1966) \cite{MOS} as
\begin{equation}
\left[\frac{\partial}{\partial\nu}P_\nu(z)\right]_{\nu=0}=\frac{z-1}{2}
{}_2F_1\left(1,1;2;\frac{1-z}{2}\right).
\label{limitderivdegreepnu0}
\end{equation}
An important generalization of this formula has recently been derived in
Szmytkowski (2011) \cite{Szmy1}. The degree-derivative of the associated
Legendre function of the first kind for $p,m\in\N_0$ and $0\le m\le p$ 
(cf.~(5.12) in Szmytkowski (2011) \cite{Szmy1}) is given by 
\begin{eqnarray}
\hspace{-0.0cm}\left[\frac{\partial}{\partial\nu}P_\nu^m(z)\right]_{\nu=p}&=&P_p^m(z)\log\frac{z+1}{2}
\nonumber\\[0.2cm]
&&\hspace{-1.5cm}{}+\left[2\psi(2p+1)-\psi(p+1)-\psi(p-m+1)\right]P_p^m(z)\nonumber\\[0.2cm]
&&\hspace{-1.5cm}{}+(-1)^{p+m}\sum_{k=0}^{p-m-1}(-1)^k\frac{(2k+2m+1)
\left[1+\frac{k!(p+m)!}{(k+2m)!(p-m)!}\right]
}{(p-m-k)(p+m+k+1)}P_{k+m}^m(z)
\nonumber\\[0.2cm]
&&\hspace{-1.5cm}{}+(-1)^{p}\frac{(p+m)!}{(p-m)!}\sum_{k=0}^{m-1}(-1)^k\frac{2k+1}{(p-k)(p+k+1)}
P_k^{-m}(z),
\label{limitderivdegreepnum}
\end{eqnarray}
and for $m\ge p+1$ (cf.~(5.16) in Szmytkowski (2011) \cite{Szmy1}) there is
\begin{equation}
\left[\frac{\partial}{\partial\nu}P_\nu^m(z)\right]_{\nu=p}=
(-1)^{p+n+1}(p+n)!(n-p-1)!P_p^{-n}(z).
\label{limitderivdegreepnum2}
\end{equation}
Some special cases of 
(\ref{limitderivdegreepnum})
include
for $m=0$
\begin{eqnarray*}
\left[\frac{\partial}{\partial\nu}P_\nu(z)\right]_{\nu=p}&=&P_p(z)\log\frac{z+1}{2}
+2\left[\psi(2p+1)-\psi(p+1)\right]P_p(z)\nonumber\\[0.2cm]
&&{}+2(-1)^{p}\sum_{k=0}^{p-1}(-1)^k\frac{2k+1}{(p-k)(p+k+1)}P_k(z),\nonumber
\end{eqnarray*}
for $m=p$
\begin{eqnarray*}
\left[\frac{\partial}{\partial\nu}P_\nu^p(z)\right]_{\nu=p}&=&P_p^p(z)\log\frac{z+1}{2}
+\left[2\psi(2p+1)-\psi(p+1)+\gamma\right]P_p^p(z)\nonumber\\[0.2cm]
&&\hspace{-1cm}{}+(-1)^{p}(2p)!\sum_{k=0}^{p-1}(-1)^k\frac{2k+1}{(p-k)(p+k+1)}
P_k^{-p}(z),
\end{eqnarray*}
where $\gamma=-\psi(1)$ is Euler's constant 
$\gamma\approx 0.57721566490$.
Of course we also have for $m=p=0$

\[
\left[\frac{\partial}{\partial\nu}P_\nu(z)\right]_{\nu=0}=\log\frac{z+1}{2},
\]
which exactly matches (\ref{limitderivdegreepnu0}).

\section{The logarithmic polynomials}
\label{Thelogarithmicpolynomials}

The logarithmic polynomials $R_p^k$ are nonvanishing only for $-p\le k\le p$
and by construction, they are satisfied by the following recurrence relation
\begin{equation}
R_p^k(x)=\frac12 R_{p-1}^{k-1}(x)+xR_{p-1}^k(x)+\frac12 R_{p-1}^{k+1}(x).
\label{recurrencerelation}
\end{equation}
\noindent 
From (\ref{logcoshcos}) we have that $R_0^0(x)=1$.
This gives us the
starting point for the recursion. 
We conjecture that the derivative of the logarithmic polynomials $R_p^k$ is given by
\[
\frac{d}{dx}R_p^{\pm k}(x)=pR_{p-1}^{\pm k}(x).
\]
It is evident by construction that these polynomials are even in the index $k$, i.e.,\
\[
R_p^k(x)=R_p^{-k}(x).
\]
\noindent Some of the first few logarithmic polynomials are given by
\begin{eqnarray}
&R_0^0(x)=1,&\nonumber\\[0.2cm]
&R_1^0(x)=x,\ R_1^{\pm 1}(x)=\frac12&\nonumber\\[0.2cm]
&R_2^0(x)=\frac12+x^2,\ R_2^{\pm 1}(x)=x,\ R_2^{\pm 2}(x)=\frac14,&\nonumber\\[0.2cm]
&R_3^0(x)=\frac32 x+x^3,\ R_3^{\pm 1}(x)=\frac38+\frac32 x^2,\ R_3^{\pm 2}(x)=\frac34 x,\
R_3^{\pm 3}(x)=\frac18.&\nonumber
\end{eqnarray}

We can find the generating function for the logarithmic polynomials as follows.  Let
\[
F(x,y,z)=\sum_{p=0}^\infty\sum_{k=-\infty}^\infty R_p^k(x)y^kz^p
\]
be the generating function for the logarithmic polynomials $R_p^k$.  If we define the function
\[
S_p(x,y)=\sum_{k=-\infty}^\infty R_p^k(x)y^k,
\]
then using the recurrence relation for $R_p^k$ 
(\ref{recurrencerelation}) we can show
\[
S_p(x,y)=
{\displaystyle \left(x+\frac12\left(y+\frac{1}{y}\right)\right)}
S_{p-1}(x,y).
\]
Combining this result along with the fact that $R_0^0(x)=1$, we have
\[
S_p(x,y)={\displaystyle \left(x+\frac12\left(y+\frac{1}{y}\right)\right)^p},
\]
so therefore the generating function for the logarithmic polynomials $R_p^k$ is given by
\[
F(x,y,z)=\frac{1}{\displaystyle 1-z\left(x+\frac12\left(y+\frac{1}{y}\right)\right)}.
\]

An algorithm for generating the logarithmic polynomials can be obtained by solving
the following set of difference equations 
\begin{equation}
\left.
\begin{array}{rcl}
a_0(p)&=&\frac12 a_0(p-1)\\[0.3cm]
a_1(p)&=&\frac12 a_1(p-1)+xa_0(p-1)\\[0.3cm]
a_2(p)&=&\frac12 a_2(p-1)+xa_1(p-1)+\frac12a_0(p-1)\\[0.3cm]
&\vdots&\\[0.3cm]
a_n(p)&=&\frac12 a_n(p-1)+xa_{n-1}(p-1)+\frac12a_{n-2}(p-1)
\end{array}
\quad\right\},
\label{algo1}
\end{equation}
subject to the boundary conditions
\begin{equation}
\left.
\begin{array}{rcl}
a_0(0)&=&1\\[0.2cm]
a_1(1)&=&xa_0(0)\\[0.2cm]
a_2(2)&=&xa_1(1)+a_0(1)\\[0.2cm]
&\vdots&\\[0.2cm]
a_n(n)&=&xa_{n-1}(n-1)+a_{n-2}(n-1)
\end{array}
\quad\right\},
\label{algo2}
\end{equation}
where $R_p^k=a_{p-|k|}(p)$ is given along diagonals for a fixed $p-|k|$.
For instance, one can obtain
$R_p^{\pm p}(x)=\frac{1}{2^p},$
for $p\ge 0$,
$R_p^{\pm(p-1)}(x)=\frac{p}{2^{p-1}}x,$
for $p\ge 1$ and
$R_p^{\pm(p-2)}(x)=\frac{p}{2^{p}}+\frac{p(p-1)}{2^{p-1}}x^2,$
for $p\ge 2,$ etc.

\section*{Acknowledgements}
I would like to thank Tom ter Elst, Matthew Auger, and Rados{\l}aw Szmytkowski
for valuable discussions.  I would also like to thank Shenghui Yang at Wolfram 
Research for valuable assistance in generating an algorithm to symbolically 
compute the logarithmic polynomials.  I acknowledge funding for time to write this 
paper from the Dean of the Faculty of Science at the University of Auckland in the 
form of a three month stipend to enhance University of Auckland 2012 PBRF Performance.
Part of this work was conducted while H.~S.~Cohl was a National Research Council
Research Postdoctoral Associate in the Information Technology Laboratory at the
National Institute of Standards and Technology, Gaithersburg, Maryland, U.S.A.




\label{lastpage}

\end{document}